\documentclass[10pt,bezier]{article}
\usepackage{amsmath,graphicx,amssymb,amsfonts}

\newtheorem{prethm}{{\bf Theorem}}

\newenvironment{thm}{\begin{prethm}{\hspace{-0.5
               em}{\bf.}}}{\end{prethm}}

\newtheorem{prepro}[prethm]{Proposition}
\newenvironment{pro}{\begin{prepro}{\hspace{-0.5
               em}{\bf.}}}{\end{prepro}}

\newtheorem{prelem}[prethm]{Lemma}

\newtheorem{precor}[prethm]{Corollary}
\newenvironment{cor}{\begin{precor}{\hspace{-0.5
               em}{\bf.}}}{\end{precor}}

\newtheorem{preremark}{Remark}
\newenvironment{remark}{\begin{preremark}{\hspace{-0.5
               em}{\bf.}}}{\end{preremark}}

\newtheorem{preexample}{{\bf Example}}

\newenvironment{example}{\begin{preexample}\em{\hspace{-0.5
               em}{\bf.}}}{\end{preexample}}

\newtheorem{preproof}{{\bf Proof.}}

\newenvironment{proof}[1]{\begin{preproof}{\rm
               #1}\hfill{$\Box$}}{\end{preproof}}

\def \bangle{ \atopwithdelims \langle \rangle}
\title{\bf\ Some remarks about the derivative operator and generalized Stirling numbers}
\author{{\normalsize{ M. Mohammad-Noori ${}^{1}$}
}\vspace{2mm} \\{\footnotesize{\it
Department of
Mathematics, Statistics and Computer
Science, University of Tehran,}}\vspace{-2mm}\\{\footnotesize{\it   Tehran,
Iran}}\vspace{-2mm}\\
{\footnotesize{\it School of Mathematics, Institute for Research in
Fundamental Sciences {\rm(IPM),}}}\vspace{-2mm}\\{\footnotesize{\it P.O.Box: 19395-5746, Tehran,
Iran}}
\\{\footnotesize Emails: morteza@ipm.ir, mnoori@khayam.ut.ac.ir}
 }
\date{}
\begin{document}
\maketitle
\footnotetext[1]{\tt This research was in part supported by a grant from IPM (No. 88050011)}
\begin{abstract}
\noindent Studying expressions of the form $(f(x)D)^p$, where $D={\displaystyle \frac{d}{dx}}$ is the derivative  operator,
 goes back to Scherk's Ph.D. thesis in 1823. We show that this can be extended as\\ {\footnotesize
${\displaystyle\sum
 \gamma_{p;a} (f^{(0)})^{a(0)+1}
(f^{(1)})^{a(1)}\cdots(f^{(p-1)})^{a(p-1)}D^{p-\sum_i i a(i)}}$}, where the summation is taken over
the $p$-tuples $(a_0, a_1, \ldots , a_{p-1})$, satisfying $\sum_{i}a(i)=p-1,\, \sum_{i}i
a(i)<p$, $f^{(i)}=D^i f$ and $\gamma_{p;a}$ is the number of increasing trees on the vertex set $[0, p]$ having
 $a(0)+1$ leaves and having $a(i)$ vertices with $i$ children for $0<i<p$. Thus, previously known results about increasing trees, lead us to some equalities containing coefficients $\gamma_{p;a}$. In the sequel, we consider the expansion of $(x^k D)^p$ and coefficients appearing there, which are called generalized Stirling numbers by physicists. Some results about these coefficients and their inverses are discussed through bijective methods.
 Particularly, we introduce and use the notion of $(p, k)$-forest in these arguments.

\vspace{3mm}

 \end{abstract}

\section{Introduction}
The derivative operator $\frac{d}{dx}$ (or briefly $D$) plays an
important role in the theory of formal power series. There are some
known results discussing about iteration of expressions of form $f(x)D$, of which, the followings are the most well-known.  Here ${p\brace m}$
(resp. ${p\brack m}$) stands for Stirling number of the second type (resp. signless Stirling
number of the first type).

\begin{align}
\label{xDp} (xD)^p &=\sum_{m=1}^p {p\brace m}\, x^m D^m, \\
\label{xpDp} x^p D^p &=\sum_{m=1}^p (-1)^{p-m} {p\brack m} (xD)^p.
\end{align}

These equations were obtained by H. Scherk in his ph.D. thesis defended in 1823 \cite{scherk}, in which,
he has discussed expressions of the form $(f(x)D)^p$. The formula (\ref{xDp}) is related to the case
$f(x)=x$, but Scherk has also considered some other special cases such as $f(x)=e^{kx}$ and $f(x)=x^k$
as well as the general case. In the general case, he has obtained exact values for the first few coefficients as well
as some infinite classes of particular coefficients.
 Recently these studies have been followed independently in quantum physics in studying creation and annihilation operators (See \cite{Blsk1,Blsk2,Lang,PhiliPBla} and the references therein). The relation between the coefficients appearing in these expansions and some combinatorial objects are also studied in these references.\\

In this paper we adopt a vector notation, which gives an easy description for the expansion of $(f(x)D)^p$
 and a simple recurrence formula for the coefficients appearing in the expansion. The relation between these coefficients and increasing trees is already known \cite{PhiliPBla}. By expanding this relation, we
 give a formula for these coefficients and we provide bijective proofs for some of identities satisfied by these coefficients; particularly, the appearance of Stirling and Eulerian numbers in these coefficients are justified.\\

We also discuss the expansion of $(x^kD)^p$ in a slightly different way from \cite{Lang}. By studying the relation between this and the operator $xD$, not only we obtain the coefficients of the expansion as polynomials in terms
of $k-1$ (with productions of Stirling numbers as their coefficients), but also we go towards a bijective proof for
this expansion. This is done using the so called $(p,k)$-forests which are in fact some ``increasing forests".
Expansion of $x^{kp}D^p$ in terms with the general form $x^{(p-m)(k-1)}(x^kD)^m$, is
also studied as in \cite{Lang} as the inverse process.

\section{Preliminaries and Notation}
In this section, we introduce some definitions and notation which is useful in the rest of the paper.
 For a positive integer $n$,
the rising factorial (resp. falling factorial), denoted as $(x)_{\overline{n}}$ (resp. denoted as $(x)_{n}$),
is defined as $(x)_{\overline{n}}=x(x+1)\ldots (x+n-1)$ (resp. defined as $(x)_n=x(x-1)\ldots(x-n+1)$). For a real number $\alpha$, by $xD+\alpha$
 we mean $xD+\alpha{\bf 1}$ where ${\bf 1}$ is the identity operator, thus the expression $(xD)_m$ is defined as $xD(xD-1)\cdots(xD-m+1)$.\\

Since the Stirling numbers and Eulerian numbers appear in these expansions, a brief
introduction to these numbers and some related formulas appear in
the sequel of this introduction, based on the notation of \cite{knuth}. We recall that the {\em signless
Stirling numbers of the first kind}, denoted by ${n\brack k}$,
are the number of permutations of $S_n$ having exactly $k$ cycles. The
{\em Stirling numbers of the first kind} are then defined as
$(-1)^{n-k}{n\brack k}$. The {\em Stirling numbers of the second
kind}, denoted by ${n\brace k}$,
are the number of ways of
partitioning an $n$-set into exactly $k$ parts.
The {\em Eulerian  numbers}, denoted by
${n \bangle k}$,
are the number of
permutations $\pi\in S_n$ which have $k$ descents
(i.e. $|\{i:1\leq i<n,\pi_i>\pi_{i+1}\}|=k$.) Both kinds of Striling numbers and the Eulerian numbers
satisfy some binomial-type recurrence relations and nice identities (See \cite{CameronBook,knuth} for instance).
 Here, the
following identity is useful.
\begin{align}
 (x)_{\overline{n}} =\sum_j {n\brack j}x^j.\label{incrFactor}
\end{align}
 It is concluded from (\ref{xpDp}) that the following equation holds
 \begin{equation}
\label{xDprev} x^p D^p=(xD)_m.
\end{equation}

The set of integers (resp. nonnegative integers) is denoted by
$\mathbb{Z}$ (resp. $\mathbb{N}$). For integers $m$ and $n$ we
denote the set $\{x\in \mathbb{Z}: m\leq x\leq n\}$ by $[m,n]$.
We denote the set of infinite row vectors of nonnegative integers by
$\mathbb{N}^{\infty}$, so each element ${\bf a}\in
\mathbb{N}^{\infty}$ is represented as ${\bf a}=({\bf a}(0),{\bf
a}(1),\cdots,{\bf a}(p-1),\cdots)$. The vectors ${\bf j},{\bf e_m},{\bf n}\in \mathbb{N}^{\infty}$ are
defined respectively by ${\bf j}(i)=1, {\bf e_m}(i)=\delta_{m,i}$
and ${\bf n}(i)=i$ for any integer $i\geq 0$. The value of $f^{\bf a}$, for a vector ${\bf a}\in \mathbb{N}^{\infty}$ with finitely many nonzero components,
is defined as
\begin{equation*}
f^{\bf a}=\prod_{i}(f^{(i)})^{{\bf a}(i)},
\end{equation*}
where $f^{(0)}=f$ and for $j\geq 1$, we have $f^{(j)}=D^j f$. Also we define the set $\Lambda_p$ by
\begin{equation*}
\Lambda_p=\{{\bf a}\in \mathbb{N}_p^{\infty}: {\bf a}.{\bf j}^{\top}=p-1,\, {\bf a}.{\bf n}^{\top}<p  \}.
\end{equation*}

Let $V$ be a finite ordered set with $v_0=\min V$ (for instance,
$V$ can be considered as a finite set of integers). An increasing
tree on $V$, is a tree $T$ rooted at $v_0$ with $V(T)=V$, such
that for any $v \in V$, the vertices in the unique $v_0-v$ path
$P$ in $T$, appear increasingly.  A {\it starlike increasing tree} is an increasing tree,in which, any vertex
(except possibly the root) has at most one child. The increasing trees
are widely studied in the literature (See Section 1.3 of
\cite{StanEnumC}; for more information see \cite{ItFlajolet}). For a vertex $v$
of an increasing tree, we denote the number of its children by $d'(v)$.\\

In Proposition \ref{alph-pm-Alternative}, a little about finite calculus and its
notations is required: If $\{a_n\}_{n\geq 0}$ is a sequence of
real (or complex) numbers, then the operator $\Delta$ on
it is defined by $\Delta a_n=a_{n+1}-a_n$.
 It is easily proved that the following identity holds:
\begin{equation*}
\Delta^p a_n=\sum_{i=0}^{p} (-1)^i {p \choose i} a_{n+p-i}
\end{equation*}
and in particular
\begin{equation*}
\Delta^p a_n\bigg]_{n=0}=\sum_{i=0}^{p} (-1)^i {p \choose i} a_{p-i}.
\end{equation*}

\section{Expanding $(f(x)D)^p$}

After testing some small cases, one can guess that $(fD)^p$ is
expressed in the following form
\begin{equation*}
 \sum_{ \sum_{i}a_{i}=p-1,\, \sum_{i}i
a_{i}<p} \gamma_{p;a_0,a_1,\cdots,a_{p-1}} (f^{(0)})^{a_0+1}
(f^{(1)})^{a_1}\cdots(f^{(p-1)})^{a_{p-1}}D^{p-\sum_i i a_i},
\end{equation*}
where the constants $\gamma_{p;a_0,a_1,\cdots,a_{p-1}}$ are nonnegative integers.
 This is formulated and proved shortly in the following theorem, by using the vector notation of Section 2.\\

\begin{thm}
\label{fDp1} Let $p$ be a positive integer.
\begin{enumerate}
  \item[\rm(i)] We have
\begin{equation}\label{fDp}
(fD)^p= \sum_{\bf a} \gamma_{p;{\bf a}}f^{{\bf a}+{\bf
e_0}}D^{p-{\bf a}.{\bf n}^{\top}}
\end{equation}

where the summation runs over the elements ${\bf a}\in
\Lambda_p$. Equivalently, one can say that ${\bf
a}$ runs over $\mathbb{N}^{\infty}$ but $\gamma_{p;{\bf a}}\neq 0$
only if ${\bf a}\in \Lambda_p$.

  \item[\rm(ii)] For any $p\geq 2$  and for any ${\bf a}\in
\mathbb{N}^{\infty}$, we have the following recurrence relation:
\begin{equation}
\label{gammarecur} \gamma_{p;{\bf a}}= \gamma_{p-1;{\bf
a-e_0}}+\sum_{i=0}^{p-2}({\bf a}(i)+1)\gamma_{p-1;{\bf
a-e_0+e_i-e_{i+1}}}, \end{equation} where by definition the terms
$\gamma_{p-1;{\bf b}}$ in the right are $0$ whenever ${\bf b}\not
\in \Lambda_{p-1}$.
  \item[\rm(iii)]  If ${\bf a}(p-1)=1$ and $\gamma_{p;{\bf a}}\neq 0$ then
  ${\bf a}=(p-2){\bf e_0}+{\bf e_{p-1}}$ and
$\gamma_{p;(p-2){\bf e_0}+{\bf e_{p-1}}}=1$. Furthermore, if ${\bf
a}(p-1)=0$ then the upper limit of the summation in
(\ref{gammarecur}) may be considered $p-3$.

\end{enumerate}

\end{thm}

\begin{proof}{

\begin{enumerate}
  \item[\rm(i)] The proof is by induction on $p$. The case $p=1$ is trivial:
$fD=(f^{(0)})^1D$ so for $p=1$ the only nonzero coefficient is
$\gamma_{1;{\bf e_0}}=1$. Now let the above equation holds for
$p-1$. To prove it for $p$, we write
\begin{align*}
(fD)^p&=fD((fD)^{p-1})\nonumber \\
&=fD\left(\sum_{\bf a} \gamma_{p-1;{\bf a}}f^{{\bf a}+{\bf
e_0}}D^{p-1-{\bf
a}.{\bf n}^{\top}}\right)\nonumber \\
&=\sum_{\bf a} \gamma_{p-1;{\bf a}} fD\left( f^{{\bf a}+{\bf
e_0}}D^{p-1-{\bf
a}.{\bf n}^{\top}}\right)\nonumber \\
&=\sum_{\bf a} \gamma_{p-1;{\bf a}} f^{{\bf a}+2{\bf
e_0}}D^{p-{\bf a}.{\bf n}^{\top}}+\\
& \sum_{\bf a} \sum_{i=0}^{p-2}\gamma_{p-1;{\bf a}}\left({\bf a}(i)+{\bf
e_0}(i)\right)f^{{\bf a}+2{\bf e_0}-{\bf e_i}+{\bf
e_{i+1}}}D^{p-1-{\bf a}.{\bf n}^{\top}}
\end{align*}

Now it is straightforward to show that all the terms in the right (ignoring the
coefficients) are of the form $f^{{\bf b}+{\bf e_0}}D^{p-1-{\bf
b}.{\bf n}^{\top}}$ for a vector ${\bf b}$ of length $p$
satisfying ${\bf b}.{\bf j}^{\top}=p-1$ and ${\bf b}.{\bf
n}^{\top}<p$.

\item[\rm(ii)] From the proof of part (i) we have {\footnotesize
\begin{align*}
& \sum_{\bf a} \gamma_{p;{\bf a}}f^{{\bf a}+{\bf e_0}}D^{p-{\bf
a}.{\bf n}^{\top}}=\sum_{\bf c} \gamma_{p-1;{\bf c}} f^{{\bf
c}+2{\bf e_0}}D^{p-{\bf c}.{\bf n}^{\top}}+\\
& \sum_{\bf c}\sum_{i=0}^{p-2}\gamma_{p-1;{\bf c}}\left({\bf c}(i)+{\bf
e_0}(i)\right)f^{{\bf c}+2{\bf e_0}-{\bf e_i}+{\bf
e_{i+1}}}D^{p-1-{\bf c}.{\bf n}^{\top}}
\end{align*}
}
 Calculating the coefficient of $f^{{\bf a}+{\bf e_0}}D^{p-{\bf
a}.{\bf n}^{\top}}$ in both sides gives the following identity
which concludes the result.
\begin{equation*}
\gamma_{p;{\bf a}}= \gamma_{p-1;{\bf a-e_0}}+\sum_{i=0}^{p-2}({\bf
a}(i)+1)\gamma_{p-1;{\bf a-e_0+e_i-e_{i+1}}}. \end{equation*}

 \item[\rm(iii)] This is easily obtained by using previous parts.

\end{enumerate}

}\end{proof}


A simple way to find a combinatorial interpretation for expression $(fD)^p$ is labeling
the terms $f$ and also $D$ from right to left and adding an element $f_0$ in the right to
obtain an expression of the form $f_p D_p f_{p-1} D_{p-1}\ldots f_1 D_1 f_0$, and then expanding
this expression and deleting the labels (except for $f_0$) at the end.
 The following theorem is then a simple observation.

\begin{pro}\label{ballsUrns}
Let $p$ be a positive integer and ${\bf a} \in
\Lambda_p$. Then the value of $\gamma_{p;{\bf a}}$
equals the number of distributions of $p$ distinguished balls $D_1,\cdots,D_p$ into $p$
distinguished urns $f_0,\cdots,f_{p-1}$
satisfying the following conditions: \\
(i) The label of each ball is greater than the label of the urn
containing it.\\
(ii) The number of the urns with positive label which contain
exactly $i$  balls is ${\bf a}(i)$, \, for $i=0,\cdots,p-1$.
\end{pro}
 The vector ${\bf a}$ is called the {\it counting vector} of the distribution.
As a simple application of the above Proposition, one can prove the recurrence relation (\ref{gammarecur}) by the combinatorial argument obtained form it. In Section 5, we use slight modification of this approach
to associate increasing trees to the coefficients $\gamma$.

\section{The operators $(e^{kx} D)^p$ and $(x^k D)^p$}

In this section we study the expansion of $(f(x)D)^p$ in two cases
$f(x)=e^{kx}$ and $f(x)=x^k$. Both cases are considered in
\cite{scherk} and \cite{Lang}, but we study the case $f(x)=x^k$ slightly differently.
 It is remarkable that Stirling numbers of both kinds appear explicitly in our formulae.

The following proposition is mentioned in p.9-10 of \cite{scherk} as well as \cite{Lang}

\begin{pro}\label{expkxDp} Let $p$ be a positive integer. Then we have
\begin{align}\label{ekxDp}
(e^{kx}D)^p&=e^{kpx} \sum_{m=1}^p {p \brack m} k^{p-m}D^m,\\
\label{ekpxDp} e^{kpx}D^p&=\sum_{m=1}^p {p \brace
m}(-k)^{p-m}e^{k(p-m)x}(e^{kx}D)^m.
\end{align}
\end{pro}

The problem of extending $(x^kD)^ p$ is a natural generalization of (\ref{xDp}) which is studied both in
\cite{scherk} and \cite{Lang}. First, we mention the following proposition which is easily proved by induction on $p$.

\begin{pro} Let $p$ be a positive integer.
\begin{enumerate}
  \item[\rm(i)] {\rm \cite[Lemma 1]{Lang}} We have
\begin{equation}
\label{zkDp1} (x^k D)^p=x^{p(k-1)}\sum_{m=1}^p \alpha_{pm}(k) x^m
D^m,
\end{equation}
where $\alpha_{pm}(k)$ satisfies the following recurrence relation
\begin{align}
\label{alphaRec1} \alpha_{pm}(k)&=\alpha_{p-1,m-1}(k)+\big((p-1)(k-1)+m\big)\alpha_{p-1,m}(k),\\
\label{alphaRecInitial} \alpha_{p,0}(k)&=0, (p>0),\,\,\,\, \alpha_{p,p}(k)=1, (p\geq 0),\,\,\,\, \alpha_{p,q}(k)=0, (0\leq p<q).
\end{align}

  \item[\rm(ii)] {\rm \cite[Lemma 14]{Lang}} We have
\begin{equation}
\label{inv} x^{kp}D^p=\sum_{m}\beta_{pm}(k) x^{(p-m)(k-1)}(x^k D)^m,
\end{equation}
where $\beta_{pm}(k)$ satisfies the following recurrence relation
\begin{align}
\label{betaRec1} \beta_{pm}(k)&=\beta_{p-1,m-1}(k)-\big(m(k-1)+p-1\big)\beta_{p-1,m}(k),\\
\label{betaRecInitial} \beta_{p,0}(k)&=0, (p>0),\,\,\,\, \beta_{p,p}(k)=1, (p\geq 0),\,\,\,\, \beta_{p,q}(k)=0, (0\leq p<q).
\end{align}
\end{enumerate}
\end{pro}

Note that from definitions of coefficients $\alpha$ and $\beta$, it is easily seen
that they are inverses of each other, i.e. $\sum_m \alpha_{pm}\beta_{mq}=\sum_m \beta_{pm}\alpha_{mq}=\delta_{pq}$.
To calculate the coefficients $\alpha_{pm}$ and $\beta_{pm}$, one can use the equations given in the above proposition,  but we prefer to do this calculation by the auxiliary tool given below.

\begin{pro}\label{xkDxDThm} Let $p$ be a positive integer. The following identity holds between the operators
$x^kD$ and $xD$.
\begin{equation*}
(x^kD)^p=x^{p(k-1)}\prod_{i=0}^{p-1}\big(xD+i(k-1)\big).
\end{equation*}
\end{pro}

\begin{proof}
{This is proved easily by induction on $p$.}
\end{proof}

In the next theorem, we calculate $\alpha_{pm}(k)$ and $\beta_{pm}(k)$ as
 polynomials of degree $p-m$ in terms of $k-1$, whose coefficients are given in
 terms of Stirling numbers as follows:

\begin{thm}\label{alphapmkValu}
 We have
\begin{align}\label{alphaFormul}
\alpha_{pm}(k)&=\sum_{j=m}^{p} {p \brack j}{j\brace m}(k-1)^{p-j}\\
 \label{beta}
\beta_{pm}(k)&=(-1)^{p-m}\sum_{j=m}^{p} {p \brack j}{j\brace
m}(k-1)^{j-m}
\end{align}
\end{thm}

\begin{proof}
{ By (\ref{zkDp1}) and Proposition \ref{xkDxDThm}, we obtain
\begin{align*}
\sum_{m=1}^p \alpha_{pm}(k) x^mD^m&=\prod_{i=0}^{p-1}\big(xD+i(k-1)\big)\\
&=\sum_{j=1}^p {p \brack j}(k-1)^{p-j} (xD)^j\\
&=\sum_{j=1}^p {p \brack j}(k-1)^{p-j} \sum_{m=1}^p {j\brace m} x^m D^m\\
&=\sum_{m=1}^p \sum_{j=1}^p  {p \brack j}{j\brace m}(k-1)^{p-j} x^m D^m.
\end{align*}

 Hence,
\begin{equation*}
\alpha_{pm}(k)=\sum_{j=1}^{p} {p \brack j}{j\brace m}(k-1)^{p-j}=\sum_{j=m}^{p} {p \brack j}{j\brace m}(k-1)^{p-j} .
\end{equation*}
Clearly (\ref{beta}) holds for $k=1$, so we prove it in the case $k\neq 1$.
From Theorem \ref{xkDxDThm}, by changing $p$ to $n$, we get
\begin{equation*} s_n=\sum_m (-1)^{n-m}{n \brack m} t_m,
\end{equation*}
where $s_n= x^{-n(k-1)}(1-k)^{-n}(x^kD)^n$ and
$t_m=(1-k)^{-m}(xD)^m$.
By inverting the above equation, we provide
 $t_n=\sum_m {n\brace m}s_m$, which equals
\begin{equation*}
(xD)^n=(1-k)^n\sum_m(1-k)^m x^{m(1-k)}{n\brace m}(x^kD)^n.
\end{equation*}
Finally combining this with (\ref{xpDp}) yields equation (\ref{beta}).
}
\end{proof}

\begin{cor}
\label{alphabeta}
{With previous notations, the following identity holds
\begin{equation*}
\beta_{pm}(k)=(1-k)^{p-m}\alpha_{pm}(\frac{k}{k-1}).
\end{equation*}
}
\end{cor}

It is immediately seen from (\ref{alphaFormul}) that the case $k=2$ gives famous Lah numbers which are defined
as $L(p,m)=\sum_j {p \brack j}{j\brace m}$, and satisfy $L(p,m)=\frac{p!}{m!}{p-1 \choose m-1}$.
(See for instance Exercise 2.13 of \cite{CameronNotes}). By using Corollary \ref{alphabeta},
$\beta_{pm}(2)=(-1)^{p-m}\alpha_{pm}(2)$. Hence
\begin{align*}
\alpha_{pm}(2)&=\frac{p!}{m!}{p-1 \choose m-1}\\
\beta_{pm}(2)&=(-1)^{p-m}\frac{p!}{m!}{p-1 \choose m-1}
\end{align*}
The following proposition may be considered as a generalization of these equations. Moreover it gives alternative formulas for computing $\alpha_{pm}(k)$ and $\beta_{pm}(k)$.

\begin{pro}\label{alph-pm-Alternative}
\begin{enumerate}
  \item[\rm(i)] Let $p\geq m$. Then for any $t$ we have
\begin{equation*}
\sum_j {p \brack j}{j\brace m}t^j=\frac{p!}{m!} n_{pmt}.
\end{equation*}
where the value $n_{pmt}$ is given by
\begin{equation*}
n_{pmt}=\sum_{\ell=0}^m (-1)^\ell {m
\choose \ell}{p+tm-t\ell-1 \choose p}
\end{equation*}
and is equal to the number of $t \times m$ matrices with nonnegative integer entries and no zero column, and with $p$ as the sum of their entries.
\item[\rm(ii)] Let ${\displaystyle A(z)=\prod_{i=0}^{p-1}\big(z+i(k-1)\big)}$. We have
\begin{equation*}
\alpha_{pm}(k)=\frac{1}{m!} \Delta^m A(z)\big]_{z=0}.
\end{equation*}
\item[\rm(iii)] Let $k\neq 1$ and ${\displaystyle B(z)=\prod_{i=0}^{p-1}\big((k-1)z+i\big)}$. Then
\begin{equation*}
\beta_{pm}(k)=\frac{1}{m!(k-1)^m} \Delta^m B(z)\big]_{z=0}.
\end{equation*}

\end{enumerate}
\end{pro}
\begin{proof}{
We prove (i); parts (ii) and (iii) are easy conclusions of this. By replacing
${j\brace m}=\frac{1}{m!}\sum_{\ell=0}^m (-1)^{\ell}{m \choose \ell} (m-\ell)^j$, we obtain
\begin{align*}
\sum_{j=1}^p {p \brack j}{j\brace m}t^j&=\sum_{j=1}^p \sum_{\ell=0}^m  \frac{1}{m!} {p \brack j} t^j(-1)^{\ell}{m \choose \ell} (m-\ell)^j\\
&=\sum_{\ell=0}^m \frac{(-1)^{\ell}}{m!} {m \choose \ell} \sum_{j=1}^p {p \brack j} t^j (m-\ell)^j\\
&=\frac{1}{m!} \sum_{\ell=0}^m (-1)^{\ell} {m \choose \ell}(t(m-\ell))_{\overline{p}}\\
&=\frac{p!}{m!} \sum_{\ell=0}^m (-1)^{\ell}{m \choose \ell}{t(m-\ell)+p-1 \choose p}\\
&=\frac{p!}{m!} n_{pmt}
\end{align*}
By using PIE, it is easily obtained that the value $n_{pmt}$ equals the number of
 $t \times m$ matrices with nonnegative integer entries and no zero column, and with
 $p$ as the sum of their entries.
}
\end{proof}

\begin{remark} {\rm From Proposition \ref{alph-pm-Alternative}, it is concluded that if $k$ is a positive integer,
then the value $\frac{m!}{p!} (k-1)^m|\beta_{pm}(k)|$ is an integer; in fact it equals the number of $(k-1)\times m$
matrices with nonnegative entries and with no zero column whose entries sum up to $p$.}
\end{remark}

\section{Increasing trees and some bijections}
 In this section, first we mention some enumerative results about increasing trees in two propositions.
 The definition of an increasing tree related to a distribution with conditions of Proposition \ref{ballsUrns},
 is then a useful connector which concludes the rest of results of this section about
 coefficients $\gamma_{p;{\bf a}}$. The definitions and notations for increasing trees are
 as mentioned in Section 2. The proposition about increasing trees is well-known and is proved by a standard
  representation of permutations as in proposition 1.3.16 of \cite{StanEnumC}.

\begin{pro}\label{enumIA}
{ Let $\mathcal{T}(p)$ be the set of increasing trees on the
vertex set $V=[0,p]$. Then
\begin{enumerate}
\item[\rm(i)] Let $\mathcal{S}_m'(p)\subseteq \mathcal{T}(p)$ be
  the set of starlike increasing trees on $[0,p]$ for which there exists a unique
  vertex $v\in [1,p]$ such that $v$ has exactly $m$ children and for any vertex $u\in [1,p]\setminus\{v\}$, $u$ is a leaf. Then
    $|\mathcal{S}_m'(p)|={p\choose {m+1}}$.

\item[\rm(ii)] If $\mathcal{S}_m(p)\subseteq \mathcal{T}(p)$ is
the set of starlike increasing trees with $m$ leaves, then
 $|\mathcal{S}_m(p)|={p\brace m}$.

  \item[\rm(iii)] If $\tau_m(p)\subseteq \mathcal{T}(p)$ is
the set of increasing trees with $m+1$ leaves, then
$|\tau_m(p)|={p \bangle m}$.
  \item[\rm(iv)] If $\mathcal{T}_m(p)\subseteq \mathcal{T}(p)$ is
  the set of increasing trees in which the root has $m$ children, then
$|\mathcal{T}_m(p)|={p\brack m}$.
\end{enumerate} }
\end{pro}

\begin{remark}\label{dfsbij} {\rm
Since the bijection used in the proof of the above proposition is required for future sections,
 we give a naive description of it (This is slightly modified with respect to the one used in \cite{StanEnumC}).
 Let $T$ b an increasing tree on $[0,p]$ rooted at $0$. The mentioned bijection is based on the sequence $\{c_i\}_{i=0}^p$ of vertices obtained from the depth first search of $T$, with the constraint that in each step, the greater child should be visited sooner.
 Conversely, having a sequence $\{c_i\}_{i=0}^p$ with $c_0=0$ and $\{c_1,\ldots,c_p\}=[1,p]$, one can recover
 the corresponding increasing tree on $[0,p]$ by constructing the array $\{$PARENT$(v)\}_{v=1}^p$ which
assigns to any vertex $v\neq 0$ in the tree its parent through the following algorithm:
\begin{enumerate}
  \item[1.] $i=0$.
  \item[2.] $i=i+1$.
  \item[3.] If $(i=1)$ then PARENT$(c_i)=0$, go to step 2.
  \item[4.] $v=c_{i-1}$.
  \item[5.] while$(c_i<v)$ let $v=$PARENT$(v)$.
  \item[6.] PARENT$(c_i)=v$.
  \item[7.] If $(i<p)$ go to step 2.
  \item[8.] Stop.
\end{enumerate}
 }
 \end{remark}

The problem of the enumeration of increasing trees on $V=[0,p]$ which satisfy
$d'(i)=\ell_i$ for $i=1,\cdots,p-1$ where $\{\ell_i\}_{i=1}^{p-1}$ is a given sequence of
nonnegative integers , is considered in the next proposition.
(Note that since $d'(p)=0$ and $d'(0)=p-\sum_{i=1}^{p-1}d'(i)$, these
values are excluded from the sequence ${\ell_i}$.)

\begin{pro}\label{iag}
{Let $\ell_1,\ell_2,\cdots,\ell_{p-1}$ be a sequence of
nonnegative integers and let $V=[0,p]$. Then

\begin{enumerate}
  \item[\rm(i)] There exists an increasing tree $T$ on $V=[0,p]$ with
  $d'_T(v)=\ell_v$ for $v=1,\cdots,p-1$ if and only if ${\displaystyle \sum_{i=j}^{p-1} \ell_i\leq
  p-j}$ for $j=1,\cdots,p-1$.
  \item[\rm(ii)] The number of increasing trees mentioned in part (i) is
  obtained as
  $\frac{g(\ell_1,\cdots,\ell_{p-1})}{\ell_1!\cdots\ell_{p-1}!}$
  where
  \begin{equation*}g(\ell_1,\cdots,\ell_{p-1})=(2-\ell_{p-1})_*
  (3-\ell_{p-1}-\ell_{p-2})_*\cdots (p-1-\sum_{i=2}^{p-1}\ell_i)_*
\end{equation*}
  and for a real number $x$, the value of $(x)_*$ is defined to be
  $x$ if $x>0$ and $0$ otherwise.
\end{enumerate} }
\end{pro}
\begin{proof}{ \begin{enumerate}
\item[\rm(i)] We denote the set of children of the vertex $i$ by
$L_i$ for $i=1,\cdots,p-1$. Obviously theses are disjoint sets and
for any $j$~, $1\leq j\leq p-1$, we have
$\bigcup_{i=j}^{p-1}L_i\subseteq [j+1,p]$. Therefore, from the
existence of an increasing tree with the given out-degree sequence, the
required inequalities hold. Conversely, if the inequalities
${\displaystyle \sum_{i=j}^{p-1} \ell_i\leq
  p-j}$ holds for $j=1,\cdots,p-1$, one can construct $L_i$'s
  as follows: The set $L_{p-1}$ is a subset of size $\ell_{p-1}$ of
  $\{p\}$ (Thus either $\ell_{p-1}=0, L_{p-1}=\emptyset$ or $\ell_{p-1}=1,
  L_{p-1}=\{p\}$.) Now suppose that $L_i$ is constructed for
  $i=p-1,p-2,\cdots,i'+1$. Then $L_{i'}\subseteq [i'+1,p]\setminus
  \bigcup_{i=i'+1}^{p-1}L_i$, but the set in the right contains
  exactly $h_{i'}=p-i'-\sum_{i=i'+1}^{p-1}\ell_i$ elements. On the other hand, from the
  given inequalities we obtain $h_{i'}\geq \ell_{i'}$ thus it is possible to construct $L_{i'}$.
 \item[\rm(ii)] It is concluded from the construction of Part (i).
\end{enumerate}
 }
\end{proof}

{\bf Definition 1. } Consider a distribution $\mathcal{D}$ of $p$ distinguishable balls
$D_1$, $D_2$, $\ldots$, $D_p$ into $p$ distinguishable urns $f_0$, $f_1$, $\ldots$, $f_{p-1}$
satisfying conditions (i) and (ii) of Proposition \ref{ballsUrns}.  We associate a graph $T(\mathcal{D})$
with the vertex set $V=[0,p]$ to the distribution $\mathcal{D}$, as as follows: If the ball $D_i$ is put into
the urn $f_j$, then $\{i,j\}$ is an edge of $T(\mathcal{D})$. It is clear that $T(\mathcal{D})$ is an increasing tree
rooted at $0$.\\

The following theorem is concluded from the above definition and Proposition \ref{ballsUrns}.

\begin{pro}\label{arbgamma}
With the conditions of Proposition \ref{ballsUrns}, the value of $\gamma_{p;{\bf
a}}$ equals the number of
increasing trees on $[0,p]$ in which  \\
(i) The number of the leaves is ${\bf a}(0)+1$.\\
(ii) The number of the nodes which have exactly $i$ children is ${\bf a}(i)$ \, for $i=1,\cdots,p$.
\end{pro}

\begin{cor} The following identities hold.
\label{recIdnt0}
\begin{enumerate}
  \item[\rm(i)] Suppose that $1\leq m\leq p-1$. Then $\gamma_{p;(p-2){\bf e_0}+{\bf e_m}}={p \choose m+1}.$
  \item[\rm(ii)] Let $1 \leq m \leq p$. Then $\gamma_{p;(m-1)
  {\bf e_0}+(p-m){\bf e_1}}={p\brace m}.$
  \item[\rm(iii)]  ${\displaystyle \sum_{{\bf a}.{\bf n}^{\top}=p-m}\gamma_{p;{\bf a}}={p\brack m}}$
  \item[\rm(iv)]  ${\displaystyle \sum_{{\bf a}.{e_0}^{\top}=m} \gamma_{p;{\bf a}}={p \bangle m}}$.
\end{enumerate}
\end{cor}

\begin{proof}
{All parts are concluded from Proposition \ref{enumIA}. We also mention that parts (i) and (ii) can
be also proved using recurrence relation (\ref{gammarecur}). Moreover, parts (ii) and (iii) may also be
concluded from expansions of $(x D)^p$ and $(e^x D)^p$.}
\end{proof}

The following theorem gives a nonrecursive formula to compute the coefficient  $\gamma_{p;{\bf a}}$.

\begin{thm} The coefficient  $\gamma_{p;{\bf a}}$ can be computed as follows
{$$\gamma_{p;{\bf a}}=\frac{1}{(0!)^{{\bf a}(0)} (1!)^{{\bf a}(1)} \ldots ((p-1)!)^{{\bf a}(p-1)}} \sum g(\ell_1, \ell_2, \ldots, \ell_{p-1})$$
Where the summation runs over all $(p-1)$-tuple $(\ell_1, \ell_2, \ldots, \ell_{p-1})$ of integers
satisfying $\{\ell_1, \ell_2, \ldots, \ell_{p-1}\}= \{{\bf a}(0) .0, {\bf a}(1).1, \ldots,{\bf a}(p-1).(p-1)\}$
(which means that the number of $i$'s appearing in the sequence $\{\ell_i\}_{i=1}^{p-1}$ is ${\bf a}(i)$
for $i=0,\cdots,p-1$).
}
\end{thm}

\begin{proof}
{The proof is straightforward by using Proposition \ref{iag} and
the definition of $\gamma_{p;{\bf a}}$. (Note that the summation
given above, contains $\frac{(p-1)!}{{\bf a}(0)!\cdots{\bf
a}(p-1)!}$ summands, some of which may equal $0$.)}
\end{proof}

\section{Expansion of $(x^k D)^p$}
The expansion of $(x^kD)^p$ is given in \cite{Lang}. Here we would like to obtain this result through a pure combinatorial discussion. An
immediate usage of the results of the previous section, suggests to bound the capacity of each urn $f_i,~(i>0)$ by $k$, but here we propose
a more useful model:
The problem of expanding $(x^kD)^p$ is related to the expansion of the following expression
$$x_{pk}x_{pk-1}\cdots x_{pk-k+1}D_p\cdots x_{ik}x_{ik-1}\cdots x_{ik-k+1}D_i \cdots x_{k}x_{k-1}\cdots x_{1}D_1,$$
which itself is related to a $(p,k)$-distribution defined as follows (Note that the rightmost position is reserved for an urn called $x_0$).\\
\\
{\bf Definition 2. } For positive integers $p$ and $k$, a {\it $(p,k)$-distribution} is a distribution of
 balls $D_i$, ($1\leq i\leq p$) into the urns $x_j$, ($0\leq j\leq pk$) such that each urn (except $x_0$ whose capacity is not bounded)
 contains at most one ball, and a ball $D_i$ can be put into an urn $x_j$ only if $\lceil\frac{j}{k}\rceil<i$.
\begin{pro} \label{alphaBall} We have
 $(x^k D)^p=x^{p(k-1)}\sum_{m=1}^p \alpha_{pm}(k) x^m
D^m,$~where $\alpha_{pm}(k)$ is the number of $(p,k)$-distributions in which the urn $x_0$ contains exactly $m$ balls.
\end{pro}
To relate a $(p,k)$-distribution to a graph, we add a new definition. Before this, we mention that
a rooted tree $T$ with root $r$ is sometimes denoted by $(T, r)$ to emphasize on the root.

{\bf Definition 3. } A {\it $(p, k)$-forest} (related to a
$(p,k)$-distribution,) is defined as graph $F$ on the
vertex set
$V(F)=\{0,1,\cdots,(p-1)k+1\}\cup\{r_1,\cdots,r_{p-1}\}$, which is
a disjoint union of starlike increasing trees $(S_0,0)$,
$(S_1,r_1)$,$\cdots$,$(S_{p-1},r_{p-1})$, where the symbols $r_i$ are
conventionally considered smaller than any positive integer, and
the following conditions hold.

\begin{enumerate}
  \item[\rm(1)] The children of $r_i$ are vertices $(i-1)k+2,(i-1)k+3,\cdots,ik$, for $i=1,\cdots,p-1$.
  \item[\rm(2)] If the ball $D_m$ is put into urn $x_0$, then $(m-1)k+1$ is a child of $0$ in $S_0$.
  \item[\rm(3)] If the ball $D_m$ is put into urn $x_n$,
  then $(m-1)k+1$ is a child of $n$.
\end{enumerate}
We mention that $(p-1)(k-1)$ edges of the forest are determined by condition (1), independent of the related  $(p,k)$-distribution. We note that the existence of vertices $r_i$ and these edges in the $(p, k)$-forest, gives it a more symmetric structure. The other edges of the forest are obtained from the related $(p,k)$-distribution.

Also the vertex $1$ is always a child of $0$.
If $1\leq u \leq (p-1)k$ and $u\not\equiv 1 (\mod k)$, then $u$ is a child of some $r_i$. Consequently, the vertices of the form $u\equiv 1 (\mod k)$, which are not in $V(S_0)$, are at distance at least $2$ from some $r_i$ in $S_i$. Furthermore, for any such vertex $u$, there exists a unique ancestor of the form $(i-1)k+r'$ with
$(2\leq r'\leq k)$.

{\bf Definition 4.} The mapping $\Omega_{pk}: F \rightarrow ((T,0),{\mathcal C}'_k)$ maps a $(p, k)$-forest $F$ into an increasing tree $(T,0)$ and an associated coloring of its edges with the color set $C=\{0,1,\cdots,k-1\}$ constructed as follows:
\begin{enumerate}
  \item[\rm(1)] For any $v>0$ such that $(v-1)k+1\in V(S_0)$, add an edge $\{0,v\}$ into $T$ and color it with color $0$.
  \item[\rm(2)] Let $i>0$. For any vertex of the form $(v-1)k+1$ in $V(S_i)$ add an edge $\{i,v\}$ in $T$. Moreover, if the vertex $(v-1)k+1$ in $S_i$ has an ancestor of the form $(i-1)k+r+1$, with $1\leq r\leq k-1$, (note that this ancestor is unique,) then color the edge $\{i, v\}$ with color $r$.
\end{enumerate}

\begin{example}
{The forest shown in Figure 1 is the corresponding forest (according to the Definition 3) to the term
 $x_{12}x_{11} x_{10} x_{9}^{\{4\}} x_{8} x_{7} x_{6} x_{5} x_{4}^{\{3\}} x_{3} x_{2}^{\{2\}} x_{1} D^{\{1\}}$
 of the expansion of $(x^3 D)^4$\\

\vspace{-0.5cm}
\begin{figure}[h]
\def\emline#1#2#3#4#5#6{%
\put(#1,#2){\special{em:moveto}}%
\put(#4,#5){\special{em:lineto}}}
\def\newpic#1{}
%
%
%
\unitlength 0.3mm
\special{em:linewidth 0.4pt}
\linethickness{0.6pt}
\begin{picture}(150,150)(-85,0)
%
\put(-4,132){\circle*{2}}
\put(-4,111){\circle*{2}}
\put(60,132){\circle*{2}}
\put(72,111){\circle*{2}}
\put(50,111){\circle*{2}}
\put(50,89){\circle*{2}}
\put(50,68){\circle*{2}}
\put(125,132){\circle*{2}}
\put(115,111){\circle*{2}}
\put(136,110){\circle*{2}}
\put(190,132){\circle*{2}}
\put(179,110){\circle*{2}}
\put(201,111){\circle*{2}}
\put(201,89){\circle*{2}}
\emline{-4}{132}{1}{-4}{111}{2}
\emline{60}{132}{1}{72}{111}{2}
\emline{60}{132}{1}{50}{111}{2}
\emline{50}{111}{1}{50}{89}{2}
\emline{50}{89}{1}{50}{68}{2}
\emline{125}{132}{1}{115}{111}{2}
\emline{125}{132}{1}{136}{110}{2}
\emline{190}{132}{1}{179}{110}{2}
\emline{190}{132}{1}{201}{111}{2}
\emline{201}{111}{1}{201}{89}{2}
\put(-5,105){\makebox(0, 0)[cc]{1}}
\put(-4,139){\makebox(0, 0)[cc]{0}}
\put(43,69){\makebox(0, 0)[cc]{7}}
\put(43,90){\makebox(0, 0)[cc]{4}}
\put(76,112){\makebox(0, 0)[cc]{3}}
\put(43,112){\makebox(0, 0)[cc]{2}}
\put(60,139){\makebox(0, 0)[cc]{$r_1$}}
\put(126,139){\makebox(0, 0)[cc]{$r_2$}}
\put(192,139){\makebox(0, 0)[cc]{$r_3$}}
\put(141,112){\makebox(0, 0)[cc]{6}}
\put(109,112){\makebox(0, 0)[cc]{5}}
\put(207,112){\makebox(0, 0)[cc]{9}}
\put(175,112){\makebox(0, 0)[cc]{8}}
\put(207,90){\makebox(0, 0)[cc]{10}}

\end{picture}
\vspace{-3.5cm}
\end{figure}
\vspace{1cm}
\begin{center}
{Figure $1$}
\end{center}

And the colored increasing tree shown in Figure 2, is the corresponding colored increasing tree (according to the Definition 4)
to the forest shown in Figure 1.

\begin{figure}[h]
\def\emline#1#2#3#4#5#6{%
\put(#1,#2){\special{em:moveto}}%
\put(#4,#5){\special{em:lineto}}}
\def\newpic#1{}
%
%
%
\unitlength 0.4mm
%
\special{em:linewidth 0.4pt}
\linethickness{0.6pt}
\vspace{-3.5 cm}
\begin{picture}(150,150)(-91,0)
%
\put(60,61){\circle*{2}}
\put(60,39){\circle*{2}}
\put(39,18){\circle*{2}}
\put(82,18){\circle*{2}}
\put(82,-3){\circle*{2}}
\put(60,39){\circle*{2}}
\emline{60}{61}{1}{60}{39}{2}
\emline{60}{39}{1}{39}{18}{2}
\emline{60}{39}{1}{82}{18}{2}
\emline{82}{18}{1}{82}{-3}{2}
\put(92,7){\makebox(0, 0)[cc]{$c(2)$}}
\put(43,31){\makebox(0, 0)[cc]{$c(1)$}}
\put(80,31){\makebox(0, 0)[cc]{$c(1)$}}
\put(70,50){\makebox(0, 0)[cc]{$c(0)$}}
\put(81,-9){\makebox(0, 0)[cc]{4}}
\put(38,12){\makebox(0, 0)[cc]{2}}
\put(87,20){\makebox(0, 0)[cc]{3}}
\put(64,41){\makebox(0, 0)[cc]{1}}
\put(60,67){\makebox(0, 0)[cc]{0}}

\end{picture}
\end{figure}
\begin{center}{Figure 2}
\end{center}

 }
\end{example}

\begin{thm}\label{bijFArb}
Using the previous notations we have
\begin{enumerate}
  \item[\rm(i)] There exists a bijection $\Xi:F\rightarrow ((T,0),\mathcal{C}'_k,(S_0,0))$.
  \item[\rm(ii)] The mapping in part (i) gives a bijective proof for the following identity
  $$\alpha_{pm}(k)=\sum_{j=m}^{p} {p\brack j}{j\brace m}(k-1)^{p-j}.$$
  \end{enumerate}
\end{thm}
\begin{proof}
{ (i) Having the $(p, k)$-forest $F$, we obtain $(S_0,0)$ as its component which contain the vertex $0$. The other components
 of the triple in the right (which are $(T,0)$ and $\mathcal{C}'_k$,) are constructed uniquely from $F$ by the mapping $\Omega_{pk}$ as mentioned before.\\
 Now suppose that $(S_0,0),$~$(T,0),$~and $\mathcal{C}'_k$ are given. To construct $F$, it is enough to construct the components $(S_i,r_i)$ for $i=1,\cdots,p-1$. For this, the following process should be followed for $i=1,\cdots p-1$: Fix $i$ and consider all of the edges $\{i,v\}$ in $T$. For any child $v$ of $i$, let $(v-1)k+1 \in V(S_i)$ and consider the vertex  $(v-1)k+1$ in $S_i$ as a descendent of $(i-1)k+r+1$, where $r$ is the color of the edge $\{i,v\}$.\\

 (ii) To calculate $\alpha_{pm}(k)$, by Proposition \ref{alphaBall}, we should calculate the number of $(p,k)$-distributions in which $x_0$ contains $m$ balls.
By part (i), this equals the number of triples $((T,0),\mathcal{C}'_k,(S_0,0))$ in which the vertex $0$ has $m$ children in $S_0$. Suppose that the vertex $0$ has $j$ children in $T$; thus $|V(S_0)\setminus\{0\}|=j$. There are then ${p\brack j}$ choices for $T$, and for any one of them, by Proposition \ref{enumIA} (ii),
there are ${j\brace m}$ choices for $S_0$. It remains to count the number of colorings $\mathcal{C}'_k$ of $T$:
There are $j$ edges of the form $\{0,v\}\in A(T)$ with color $0$; For any of the remaining $p-j$ edges we have $k-1$ choices. Thus there are $(k-1)^{p-j}$ such colorings.
We conclude that there are ${p\brack j}{j\brace m}(k-1)^{p-j}$ possibilities for such a triple.
Finally, summing up over $j$, gives the desired value.
}
\end{proof}

{\bf Acknowledgement.}  I thank Narges Ghareghani very much for several useful ideas, suggestions and comments. Also I would like to appreciate Professor Philippe Flajolet so much for many helpful comments and sending me references \cite{PhiliPBla,Blsk1,Lang,Blsk2}.

\end{document}